\documentstyle[bezier]{article} 

\newcommand{\SL}{\mbox{$\cal S\!L$}}
\newcommand{\C}{\mbox{$\cal C\!L$}}

\newcommand{\DD}{\mbox{$|D_{2}(G)|$}}

\begin{document}
{\footnotesize This paper is published in Bulletin of the ICA, Vol. 75 (2015), 47-63. However, some lines in Fig 4.1 and fraction expression in the printout in the Bulletin are missing}
 
\title{Properties of Catlin's reduced graphs and supereulerian graphs}

\author{
Wei-Guo, Chen, Guangdong Economic Information Center\\  Guangzhou, P. R. China\\
Zhi-Hong Chen\thanks{Email: chen@butler.edu}, Butler University\\ Indianapolis, IN 46208, USA. \\
Mei Lu, Tsinghua University\\ Beijing, P. R. China
}

\date{}

\maketitle
\begin{abstract}
A graph $G$ is called collapsible if for every even subset $R\subseteq V(G)$, 
there is a spanning connected subgraph $H$ of $G$ such that
$R$ is the set of vertices of odd degree in $H$. A graph is the reduction of $G$ 
if it is obtained from $G$ by contracting all the nontrivial
collapsible subgraphs. A graph is reduced if it has no nontrivial collapsible subgraphs.
In this paper, we first prove a few results on the properties of reduced graphs. As an application,  
 for  3-edge-connected graphs $G$ of order $n$ with
 $d(u)+d(v)\ge 2(n/p-1)$ for any $uv\in E(G)$ where $p>0$ are given, we show how such graphs change if
they have no spanning Eulerian subgraphs when $p$ is increased from $p=1$ to 10 then to $15$.
\end{abstract}

\noindent {\bf 1. Introduction}
\vskip 0.1 in
We shall use the notation of Bondy and Murty \cite{BonMur}, except when otherwise stated.
Graphs considered in this paper are finite and loopless, but multiple edges are allowed.
The graph of order 2   and size 2 is called a 2-cycle and denoted by $C_2$.
As in \cite{BonMur}, $\kappa'(G)$
and $d_G(v)$ (or $d(v)$)  denote  the edge-connectivity of $G$ and the degree of a vertex $v$ in $G$, respectively.
 The size of a maximum matching in $G$ 
is denoted by $\alpha'(G)$. 
A connected graph $G$ is {\it Eulerian} if the degree of each vertex in $G$ is even. 
An Eulerian subgraph $H$ of $G$  is called a {\it spanning Eulerian subgraph} if $V(G)=V(H)$ 
and is called a {\it dominating Eulerian subgraph} if $E(G-V(H))=\emptyset$.
A graph is {\it supereulerian} if it contains a spanning Eulerian subgraph. 
The family of supereulerian graphs is denoted by $\SL$.
\vskip 0.1 in

Let $O(G)$ be the set of vertices of odd degree in $G$.
A graph $G$ is {\it collapsible} if for every even subset $R\subseteq V(G)$, there is a spanning connected subgraph $H_R$ of $G$ with 
$O(H_R)=R$.   $K_{3,3}-e$ and $K_n$  ($n\ge 3$) are collapsible \cite{Ca1}.
$K_1$ is regarded as  collapsible and  supereulerian, and having $\kappa'(K_1)=\infty$.
 The family of collapsible graphs is denoted by $\cal C\!L$ . 
Thus, ${\C\subset \SL}$.

\vskip 0.1 in

Throughout this paper, we use $P$ for the Petersen graph and use  $P_{14}$ and $P_{16}$ for the graphs  defined in Figure 1.1. 

\vskip 0.1cm

$$
\setlength{\unitlength}{.75pt}
\begin{picture}(300,80)(0,-10)
\put(-35,50){$P_{14}$}
\multiput(0,50)(14,0){2}{\circle*{4}}
\put(66,50){\circle*{4}}
\put(81,50){\circle*{4}}

\multiput(0,50)(14,0){2}{\line(1,0){20}}
\put(40,80){\circle*{4}}
\put(14,15){\circle*{4}}
\multiput(26,30)(28,0){2}{\circle*{4}}
\put(81,30){\circle*{4}}
\put(81,70){\circle*4}

\put(101,40){\circle*4}
\put(101,60){\circle*4}
\put(101,40){\line(-2,-1){20}}
\put(101,40){\line(-2,1){20}}
\put(101,40){\line(-2,3){20}}

\put(101,60){\line(-2,-1){20}}
\put(101,60){\line(-2,1){20}}
\put(101,60){\line(-2,-3){20}}

\put(66,15){\line(1,1){15}}
\put(66,50){\line(1,0){15}}
\put(0,50){\line(2,-5){14}}
\put(0,50){\line(4,3){40}}
\put(40,65){\circle*{4}}
\put(40,65){\line(0,1){15}}
\put(26,30){\line(2,5){14}}
\put(40,65){\line(2,-5){14}}
\put(26,30){\line(2,1){40}}
\put(14,50){\line(2,-1){40}}
\put(14,50){\line(1,0){52}}

\put(66,15){\circle*{4}}

\put(54,30){\line(4,-5){12}}
\put(14,15){\line(3,4){11}}
\put(14,15){\line(1,0){52}}
 
\put(40,80){\line(4,-1){40}}


\setlength{\unitlength}{1pt}

\put(180,40){\circle*3}

\put(180,40){\line(1,-3){10}}
\put(180,40){\line(1,0){10}}
\put(180,40){\line(3,2){30}}

\put(190,40){\circle*3}
\put(190,40){\line(1,0){40}}
\put(190,40){\line(3,-2){30}}

\put(230,40){\circle*3}

\put(190,10){\circle*3}
\put(190,10){\line(1,1){10}}
\put(190,10){\line(1,0){40}}

\put(230,10){\circle*3}
\put(230,10){\line(-1,1){10}}
\put(230,10){\line(1,0){20}}

\put(200,20){\circle*3}
\put(200,20){\line(1,3){10}}
\put(200,20){\line(3,2){30}}

\put(220,20){\circle*3}

\put(210,50){\circle*3}
\put(210,50){\line(0,1){10}}
\put(210,50){\line(1,-3){10}}

\put(210,60){\circle*3}
\put(250,10){\circle*3}   
\put(250,10){\line(0,1){20}}
\put(250,10){\line(2,3){20}}

\put(250,30){\circle*3}

\put(240,40){\circle*3}
\put(240,40){\line(-1,0){10}}
\put(240,40){\line(1,1){10}}
\put(240,40){\line(1,-1){10}}

\put(250,50){\circle*3}

\put(250,60){\circle*3}

\put(250,60){\line(-1,0){40}}
\put(250,60){\line(0,-1){10}}
\put(250,60){\line(1,-1){20}}

\put(260,40){\circle*3}

\put(260,40){\line(-1,1){10}}
\put(260,40){\line(-1,-1){10}}
\put(260,40){\line(1,0){10}}

\put(270,40){\circle*3}

\put(150,35){ $P_{16}$}

\put(20, -5){ (a) }
\put(210, -5){ (b)}
\put(90, -12) {Figure 1.1}
\end{picture}
$$

\vskip 0.1 in

Like the study of many NP-complete problems in graph theory,
various degree conditions for the existence of spanning and dominating Eulerian subgraphs in graphs have been derived
(e.g, see \cite{BCKV86, Brualdi, Ca1, Cat2, chen1, c4,LLZ, LaiH98, V1}).
For a graph $G$, we define\\
 \indent $\delta(G)=\min\{d(v)\ | v\in V(G)\}$;\\ 
\indent $\sigma_2(G)=\min\{d(u)+d(v)\ | uv\not \in E(G)\}$; \\ 
\indent $\sigma_{t}(G)=\min\{ \Sigma_{i=1}^{t} d(v_i)\ | \ \{ v_1, v_2, \cdots, v_t\}\ \mbox{is independent in $G$}\  (t\ge 2) \ \}$;\\ 
\indent $\delta_F(G)=\min\{ \max\{ d(u),\ d(v) \}\ |\ \mbox{ for any $u, v \in V(G)$ with $dist(u, v)=2$}\}$;\\ 
\indent $\overline{\sigma}_2(G)=\min\{ d(u)+d(v) \ | \ \mbox{for every edge $uv\in E(G)$} \}$;\\ 
\indent $\delta_L(G)=\min\{ \max \{d_G(u), d_G(v)\} | \ \mbox{ for every edge } uv \in E(G)\}$. 
\vskip 0.1 in
These are all the degree parameters we know that have been studied by many for problems on spanning and dominating Eulerian subgraphs
in graphs. In the following, we let 
$$\Omega(G)=\{\delta(G), \sigma_2(G), \sigma_{t}(G), \delta_F(G), \overline{\sigma}_2(G), \delta_L(G)\}.$$

\vskip 0.1 in

A powerful tool to work on spanning and dominating Eulerian subgraphs is Catlin's reduction method \cite{Ca1}. 
This reduction method has been applied to solve problems in Hamiltonian cycles in claw-free graphs \cite{LSZ}, 
 hamiltonian line graphs, a certain type of double cycle cover \cite{Ca4}
and the total interval number of a graph \cite{Ca5}, and others \cite{CaS}. 
\vskip 0.1 in 

\noindent {\bf Catlin's reduction method}

For $X\subseteq E(G)$, the {\it contraction} $G/X$ is the graph obtained from $G$
by identifying the two ends of each edge $e\in X$ and deleting the resulting loops. If $H$ is a subgraph of $G$, then we write $G/H$ for 
$G/E(H)$ and use $v_H$ for the vertex in $G/H$ to which $H$ is contracted.
A contraction $G/H$ is called a trivial contraction if $H=K_1$. 

\vskip 0.1 in
 Catlin \cite{Ca1} showed that every graph $G$ has a unique collection of pairwise disjoint 
maximal collapsible subgraphs $H_1$, $H_2$, $\cdots$, $H_c$ such that $V(G)=\cup_{i=1}^cV(H_i)$. 
The contraction of $G$ obtained from $G$ by contracting each $H_i$ into a single vertex 
 $v_i$ ($1\le i\le c$) is called the {\it reduction} of $G$ and denoted by $G'$. For a vertex $v\in V(G')$, 
there is a unique maximal collapsible subgraph in $G$, denoted by $H(v)$, such that $v$ is the contraction image of $H(v)$. We call 
$H(v)$  the {\it preimage} of $v$. 
A graph $G$ is {\it reduced} if $G=G'$.  By the definition of contraction, we have $\kappa'(G')\ge \kappa'(G)$. 
If the reduction of a graph $G_A$ is a graph $G_B$, we said that graph $G_A$ can be reduced to graph $G_B$.

\vskip 0.1 in
The main theorem of Catlin's reduction method is the following:\\
\noindent {\bf Theorem A} (Catlin \cite{Ca1}). Let $G$ be a graph, and let $G'$ be the reduction of $G$. Let $H$ be a collapsible subgraph 
of $G$. Then each of the following holds:
\vspace{-.1in} 
\begin{enumerate}
  \item[(a)] $G\in \C$ if and only if $G/H\in \C$. In particular, $G\in \C$ if and only if $G'=K_1$.
\vspace{-.1in} 
 \item[(b)] $G\in \SL$ if and only if $G/H\in \SL$. In particular, $G\in \SL$ if and only if $G'\in \SL$.
\end{enumerate}

\vskip 0.1 in
With Theorem A, we can see that  to determine if a graph 
 is supereulerian can be reduced to a problem of 
the reduction of the graph. For instance, 
 by combining the prior results in \cite{Cat2, chen1, c4, CL2} and the results proved recently in \cite{chenFan, chenLai}, we have:\\
\noindent {\bf Theorem B}. Let $G$ be a 3-edge-connected graph of order $n$. Let $p>1$ and $\epsilon$ be given numbers.
 Let $D(G)\in \Omega(G)$. 
If $D(G)\ge \frac{n}{p}-\epsilon$, 
then when $n$ is large,
either $G\in \SL$ or $G'$ has order at most $cp$ where $c$ is a constant.

\vskip 0.1 in

To be more specific, let $D(G)=\delta_F(G)$, we have\\
\noindent {\bf Theorem C} (W. Chen and Z. Chen \cite{chenFan}). Let $G$ be a 3-edge-connected graph of order $n$ with girth $g \in \{3,4\}$. Let
$G'$ be the reduction of $G$. 
If $\delta_F(G)> \frac{n}{(g-2)p}-\epsilon$
where $p\ge 2$ and $\epsilon>0$ are fixed and $n$ is large,
then either $G\in\SL$ or $G'\not=K_1$  has order at most $5(p-2)$. 

\vskip 0.1 in
For $D(G)=\overline{\sigma}_2(G)$, we have \\
\noindent {\bf Theorem D} (Chen and Lai \cite{chen1, CL2}). 
Let $p>0$ be an integer.
Let $G$ be a 3-edge-connected simple graph of order $n$.
 Let $G'$ be the reduction of $G$.
If $n\ge 12p(p-1)$ and 
$\overline{\sigma}_2(G)\ge \frac{2n}{p}-2,$
then either $G\in \SL$ or $G'\not =K_1$  with $\alpha'(G')\le p/2$ and $|V(G')|\le 3p/2-4$. $\Box$

\vskip 0.1 in

With Theorems B, C and D, the problem to determine if a  graph $G$ with $D(G)\ge \frac{n}{p}-\epsilon$ is in $\SL$ 
 can be reduced to
the problem of a finite number of reduced graphs. The main challenge to solve such problems  
become  solving the problems of  reduced graphs.

\vskip 0.1 in 

In this paper, we first prove some results on the properties and structures of reduced graphs. Then as an application,
we  prove a result on $\overline{\sigma}_2(G)\ge \frac{2n}{p}-2$  conditions for 3-edge-connected graphs.  
 Combining  prior results on  $\overline{\sigma}_2(G)$ conditions, it reveals how such graphs are change 
from  supereulerian to graphs that can be reduced to the Petersen graph and then to graphs that can be reduced to $P_{14}$ 
when $p$ is increased from $1$ to $10$ then to $15$.

\vskip 0.1 in

\noindent {\bf 2. Prior theorems on Catlin's reduction and $\pi$-reduction methods}

\vskip 0.1 in

For a graph $G$, let $F(G)$ be the minimum number of extra edges that must be added to $G$, 
to obtain a spanning supergraph having two edge-disjoint spanning trees. 

\noindent {\bf Theorem E}. Let $G$ be a connected reduced graph. Then each of the following holds:
\vspace{-.1in}
\begin{enumerate}
  \item[(a)] \cite{Ca1} $G$ is simple and $K_3$-free with $\delta(G)\le 3$. 
       Any subgraph $H$ of $G$ is reduced. 
\vspace{-.1in}
\item[(b)] \cite{Ca2} $F(G)=2|V(G)|-|E(G)|-2$.
\vspace{-.1in}
\item[(c)] \cite{CaHan} If $F(G)\le 2$, then $G\in \{K_1, K_2, K_{2,t} (t\ge 1)\}$. 
\vspace{-.1in}
\item[(d)] \cite{CL2} If $\delta(G)\ge 3$, then $\alpha'(G)\ge (|V(G)|+4)/3$. 
\end{enumerate}
\vskip 0.1 in

For a graph $G$, define $D_i(G)=\{ v\in V(G)\ | \ d(v)=i\}$.\\
\noindent {\bf Theorem F} (Chen \cite{c2, chen}). Let $G$ be a connected simple graph of order $n$ with $\delta(G)\ge 2$. 
Let $G'$ be the reduction of $G$. Then
each of the following holds:
\begin{enumerate}
\vspace{-.1in}
\item[(a)] \cite{c2} If $n\le 7$, $\delta(G)\ge 2$ and $\DD \le 2$, then $G$ is not reduced and   $G' \in \{ K_1, K_2 \}.$ 
\vspace{-.1in}
\item[(b)] \cite{chen} If $\kappa'(G)\ge 3$ and  $n\le 14$, then either $G\in \SL$  or  $G'\in \{P, P_{14}\}$.
\vspace{-.1in}
\item[(c)]  \cite{chen} If $\kappa'(G)\ge 3$, $n=15$, $G\not \in \SL$ and $G'\not\in \{P, P_{14}\}$, 
   then $G=G'$ has girth at least 5 and  $V(G)=D_3(G)\cup D_4(G)$ where 
$D_4(G)$ is an independent set with $|D_4(G)|=3$. 
\end{enumerate}

\vskip 0.1 in

\noindent {\bf Catlin's $\pi$-reduction method} \cite{Ca2}:
Let $G$ be a graph containing an induced 4-cycle $uvzwu$ and let $E=\{uv, vz, zw, wu\}.$
Denote by  $G/\pi$ the graph obtained from $G-E$ by identifying $u$ and $z$ to form a vertex $x$, 
and by identifying $v$ and $w$ to form a vertex $y$, and by adding an edge $e_{\pi}=xy$. 
The way to obtain $G/\pi$ from $G$ is called {\it $\pi$-reduction method} (Catlin \cite{Ca2}). 
\vskip 0.1 in

\noindent {\bf Theorem G} (Catlin \cite{Ca2}). 
Let $G$ be a connected graph and let $G/\pi$ be the graph defined above, then each of the following holds:\\
  \indent (a) If $G/\pi\in \C$, then $G\in \C$;\\
  \indent (b) If $G/\pi\in \SL$ then $G\in \SL$. $\Box$

\vskip 0.1 in

$$
\begin{picture}(380,70)(0,0)
\put(170,0){Figure 2.2}
\multiput(80,50)(0,10){3}{\circle*3}
\multiput(100,55)(0,10){2}{\circle*3}
\put(100,55){\line(-4,1){20}}
\put(100,55){\line(-4,-1){20}}
\put(100,55){\line(-4,3){20}}

\put(100,65){\line(-4,1){20}}
\put(100,65){\line(-4,-1){20}}
\put(100,65){\line(-4,-3){20}}

\put(100,65){\line(1,0){20}}    
\put(100,65){\line(2,-1){20}}
\put(100,55){\line(1,0){20}}
\put(100,55){\line(2,1){20}}
\put(98,67){\footnotesize $u$}
\put(98,50){\footnotesize $z$}

\put(120,65){\circle*3}   
\put(120,65){\line(2,1){12}}
\put(120,55){\circle*3}   
\put(120,55){\line(2,-1){12}}

\put(118, 67){\footnotesize $v$}
\put(115, 49){\footnotesize $w$}

\put(60,15){\dashbox(86,75)}
\put(90,25){\small in $G$}
\put(180,57){$\Longleftrightarrow$}


\put(255,80){\footnotesize $\Phi(x, 3)$}
\put(265,78){\vector(0,-1){10}}

\multiput(250,50)(0,10){3}{\circle*3}

\put(270,60){\circle*3}
\put(270,62){\footnotesize $x$}
\put(270,60){\line(1,0){20}}
\put(277,55){\footnotesize $e_{\pi}$}

\put(290,60){\circle*3}
\put(290,60){\line(2,1){12}}
\put(290,60){\line(2,-1){12}} 
\put(289,64){\footnotesize $y$}

\bezier{360}(250,50)(258,58)(270,60)
\bezier{360}(250,50)(258,48)(270,60)

\bezier{360}(250,60)(258,63)(270,60)
\bezier{360}(250,60)(258,57)(270,60)

\bezier{360}(250,70)(260,73)(270,60)
\bezier{360}(250,70)(258,62)(270,60)

\put(230,15){\dashbox(86,75)}
\put(250,25){\small in $G/\pi$}

\end{picture}
$$

Let $\Phi(v, t)$ be the graph obtained from $K_{1,t}$ with center at $v$ by replacing each edge in $K_{1,t}$ by a $C_2$. 
Thus, $\Phi(v,t)$ is a graph formed by $t$ $C_2$s with all the edges incident with $v$ and $|V(\Phi(v,t))|=t+1$ and $|E(\Phi(v,t))|=2t$. 
(See $\Phi(x, 3)$  in Figure 2.2). 

\vskip 0.1 in

\noindent {\bf Lemma 2.1}. Let $G$ be a connected reduced graph with $\delta(G)\ge 3$. Let $H=uvzwu$ be a 4-cycle in $G$. Let
$G/\pi$ be the graph defined by $\pi$-reduction on $G$ with $e_{\pi}=xy$. Then $G/\pi$ has at most two nontrivial collapsible subgraphs.
Furthermore, if $H_0$ is a nontrivial maximal collapsible subgraph of $G/\pi$, then 
$|V(H_0)\cap \{x,y\}| =1$ and either $H_0 =\Phi(v, t)$ 
for some $t\ge 1$ ($v\in \{x, y\}$) and $2|V(H_0)|-|E(H_0)|=2$,  or $3\le 2|V(H_0)|-|E(H_0)|$. Hence,
\begin{eqnarray*}
2\le 2|V(H_0)|-|E(H_0)|.
\end{eqnarray*}
\noindent {\bf Proof}.  Since $G$ is reduced with $\delta(G)\ge 3$, by Theorem A and Theorem G,  
$G\not =K_1$ and $(G/\pi)'\not =K_1$. If $G/\pi$ is not reduced, let $H_0$ be a nontrivial maximal collapsible 
subgraph of $G/\pi$. If $V(H_0)\cap \{x, y\}=\emptyset$, then $H_0$ is a nontrivial collapsible subgraph of $G$, contrary
to that $G$ is reduced. If $\{x, y\}\subseteq V(H_0)$, then by Theorem G, $G[E(H)\cup \{uv, vz, zw, wu\}]$ is a nontrivial collapsible subgraph 
of $G$, a contradiction again. Thus, any nontrivial maximal 
collapsible subgraph of $G/\pi$ must contain one and only one vertex in $\{x, y\}$. 

We may assume  $x\in V(H_0)$. Then $G$ has a subgraph $H_1$ with $V(H_1)=(V(H_0)-\{x\})\cup\{u, z\}$ and $E(H_1)=E(H_0)$.
If $H_0\not =\Phi(x, t)$ ($t\ge 1$), then $H_1\not =K_{2, t}$. 
Since $H_0$ is nontrivial, $H_1\not =K_2$. By Theorem E(c), $F(H_1)\ge 3$. 
Then 
\begin{eqnarray*}
3\le F(H_1)&=&2|V(H_1)|-|E(H_1)|-2\\
 &=&2(|V(H_0)|+1)-|E(H_0)|-2=2|V(H_0)|-|E(H_0)|.
\end{eqnarray*}  
Lemma 2.1 is proved. $\Box$

\vskip 0.1 in

\noindent {\bf 3. Properties of Catlin's reduced graphs}
\vskip 0.1 in

Catlin had the following conjectures on  reduced graphs:\\
\noindent {\bf Conjecture A} (Conjecture 4 of \cite{Ca4}). 
A 3-edge-connected nontrivial reduced graph $G$ with $F(G)=3$ must be the Petersen graph $P$.\\
\noindent {\bf Conjecture B} (\cite{Ca5}).
A 3-edge-connected simple graph $G$ of order at most 17 is either in $\SL$ or its reduction is in 
$\{P, P_{14}, P_{16}\}$. Thus, either $G\in \SL$ or $G$ can be contracted to $P$.

 Theorem F(b)  indicates that these conjectures are valid
 for graphs with at most 14 vertices. In this section, we prove some results 
on certain structure properties of reduced graphs that are related to these 
 conjectures and that will be needed in section 4.

\vskip 0.1 in
For convenience, for a connected graph $G$, we define 
$$f(G)=2|V(G)|-|E(G)|-2.$$
By Theorem E(b), if $G$ is reduced, then $F(G)=f(G)$. 

\vskip 0.1 in

\noindent {\bf Theorem 3.1}. Let $G$ be a connected reduced graph with $F(G)=3$ and $\delta(G)\ge 3$. 
If $G\not \in \SL$, then $G$ has no 4-cycles.\\
{\bf Proof}. 
By way of contradiction, suppose that $G$ has a 4-cycle $H_0=uvzwu$. Using $\pi$-reduction method, we have $G/\pi$ from $G$ with
$e_{\pi}=xy$ and 
\begin{eqnarray}
|V(G/\pi)|&=&|V(G)|-2\ {\rm and}\  |E(G/\pi)|=|E(G)|-3.
\end{eqnarray}
By Theorem G and the definition of $G/\pi$, since $G\not \in \SL$ with $\delta(G)\ge 3$, $G/\pi\not \in \SL$ with  $\delta(G/\pi)\ge 3$.
 By (1) and $F(G)=3$, 
\begin{eqnarray}
f(G/\pi)&=&2|V(G/\pi)|-|E(G/\pi)|-2 \nonumber \\
        &=&2(|V(G)|-2)-(|E(G)|-3)-2 \nonumber \\
        &=&2|V(G)|-|E(G)|-2-1=F(G)-1=2.
\end{eqnarray}

\vskip 0.1 in

\indent If $G/\pi$ is reduced, then by Theorem E(b) and  (2),
$F(G/\pi)=f(G/\pi)=2.$
By Theorem E(c),  $G/\pi\in \{K_1, K_2, K_{2,t}\}$, contrary to $\delta(G/\pi)\ge 3$. 
Thus, $G/\pi$ is not reduced. 

By Lemma 2.1, we may assume $G/\pi$ has a maximal collapsible subgraph $H_x$ with $x\in V(H_x)$.  By Lemma 2.1,
\begin{eqnarray}
2\le 2|V(H_x)|-|E(H_x)|.
\end{eqnarray}
Let $G_x=(G/\pi)/H_x$. Since $G/\pi\not \in \SL$, by Theorem A, $G_x\not =K_1$. 
Let $v_x$ be the vertex in $G_x$ obtained from $G/\pi$ by contracting $H_x$.
 Since $\delta(G/\pi)\ge 3$, all the vertices in $G_x$ have degree at least 3 except  $v_x$
as the result of contracting $H_x=C_2$.
By (2) and (3),
\begin{eqnarray*}
f(G_x)&=&2|V(G_x)|-|E(G_x)|-2\\
      &=&2(|V(G/\pi)|-|V(H_x)|+1)-(|E(G/\pi)|-|E(H_x)|)-2\\
      &=&f(G/\pi)+2-(2|V(H_x)|-|E(H_x)|)\le f(G/\pi)=2.
\end{eqnarray*}
If $G_x$ is reduced, then by Theorem E(c) $G_x\in \{K_1,K_2, K_{2,t}\}$, contrary to that  
 all the vertices in $G_x$ except at most one vertex have degree at least 3. 
Then $G_x$ cannot be reduced. 

 Let $H_y$ be the another nontrivial maximal collapsible subgraph of $G/\pi$. 
By Lemma 2.1, $G/\pi$ has at most two nontrivial maximal collapsible subgraphs.
Then $G_{xy}=G_x/H_y=((G/\pi)/H_x)/H_y$ is reduced.
Similar to the way of finding $f(G_x)\le 2$, we have $f(G_{xy})\le f(G_x)\le 2$ and so $F(G_{xy})=f(G_{xy})\le 2$. 
By Theorem E(c), $G_{xy}\in \{K_1, K_2, K_{2,t}\}$ ($t\ge 1$). 

If $G_{xy}=K_1$, then by Theorem A, $G/\pi\in \C\subseteq \SL$, contrary to $G/\pi\not \in \SL$.

If $G_{xy}=K_2$, then $G$ has two subgraphs $H_1$ and $H_2$ such that 
$\{u,z\}\subseteq V(H_1)$ and $E(H_1)=E(H_x)$ and $V(H_1)=(V(H_x)-\{x\})\cup \{u,z\}$, and
  $\{v, w\}\subseteq V(H_2)$ and $E(H_2)=E(H_y)$ and $V(H_2)=(V(H_y)-\{y\})\cup\{v,w \}$.
Therefore, $|E(G)|=|E(H_1)|+|E(H_2)|+|E(H_0)|=|E(H_1)|+|E(H_2)|+4$ and $|V(G)|=|V(H_1)|+|V(H_2)|$. Then
\begin{eqnarray*}
F(H_1)+F(H_2)&=&(2|V(H_1)|-|E(H_1)|-2)+(2|V(H_2)|-|E(H_2)|-2)\\
                      &=&2(|V(H_1)|+|V(H_2)|)-(|E(H_1)|+|E(H_2)|+4)\\
                      &=&(2|V(G)|-|E(G)|-2)+2=F(G)+2=5.
\end{eqnarray*}
We may assume $F(H_1)\le 2$.  Since $H_1$ is reduced, by Theorem E(c), $H_1\in \{K_1, K_2, K_{2,t}\}$. Since $H_x$ is a 
nontrivial maximal collapsible subgraph in $G/\pi$ and $G$ is reduced, $H_1\not\in \{K_1, K_2\}$. 
Hence $H_1=K_{2,t}$. Then $H_1$ has a degree two vertex $v_0\not\in \{u, z\}$. Then $d_H(v_0)=d_G(v_0)=2$, contrary 
to $\delta(G)\ge 3$. Thus, $G_{xy}=K_2$ is impossible.

If $G_{xy}=K_{2,t}$, then since $\delta(G/\pi)\ge 3$ and $K_{2,t}$ ($t\ge 1$) has at least 3 vertices with degree less than 3,
 $G/\pi$ has at least 3 nontrivial maximal collapsible subgraphs, a contradiction. Theorem 3.1 is proved. $\Box$

\vskip 0.1 in

\noindent {\bf Lemma 3.2}. Let $G$ be a connected reduced graph of order $n$. Let $H$ be a spanning bipartite subgraph of $G$
 with bipartition $\{X, Y\}$ where $|Y|\ge |X|$ and $d_H(v)\ge 3$ for any $v\in Y$.
If $|X|\le \frac{n+5}{3}$, then $G=H$ and $F(G)=3$.\\
\noindent {\bf Proof}. Since $|Y|\ge |X|$ and $d_H(v)\ge 3$ for any $v\in Y$,
$|E(H)|\ge 3|Y|$ and $|X|\ge 3$. Hence  $H\not \in \{K_1, K_2, K_{2,t}\}$ and 
so $G\not \in \{K_1, K_2, K_{2,t}\}$  $(t\ge 1)$ . 
By Theorem E(c), $F(G)\ge 3$.
Since $E(H)$, $E(G[X])$ and $E(G[Y])$ are pairwise disjoint subsets of $E(G)$, 
\begin{eqnarray}
|E(G)|&\ge& |E(H)|+|E(G[X])|+|E(G[Y])|\nonumber \\
      &\ge& 3|Y|+|E(G[X])|+|E(G[Y])|.
\end{eqnarray} 
By Theorem E(b),  (4), $|Y|\le n-|X|$ and $|X|\le \frac{n+5}{3}$,
\begin{eqnarray*}
3\le F(G)&=&2|V(G)|-|E(G)|-2\\
         &\le& 2(|X|+|Y|)-3|Y|-(|E(G[X])|+|E(G[Y])|)-2\\
         &=& 3|X|-n-2-(|E(G[X])|+|E(G[Y])|)\\
         &\le& 3(\frac{n+5}{3})-n-2-(|E(G[X])|+|E(G[Y])|)\\
         &=&3-(|E(G[X])|+|E(G[Y])|).
\end{eqnarray*}
Thus, $|E(G[X])|+|E(G[Y])|=0$ and so $G=H$ and $F(G)=3$.  Lemma 3.2 is proved. $\Box$

\vskip 0.1 in

Several properties on reduced bipartite graphs are given in the following.

\noindent {\bf Theorem 3.3}. Let $G$ be a 3-edge-connected reduced graph.
 Let $H$ be a connected reduced bipartite graph with bipartition 
$\{X, Y\}$ where $|X|\le 7$, $|Y|\ge |X|$ and $d_H(v)\ge 3$ for any $v\in Y$.
\begin{enumerate}
  \item[(a)] If $|Y|\ge |X|$, then either $H$ has a 4-cycle with a vertex of degree at least 4 in $X$ or $|Y|=|X|$ and $H\in \SL$.
  \item[(b)] If $|Y|=|X|$ and $|X|\le 6$, then $H$ has a 4-cycle.
  \item[(c)] If $G$ has such a bipartite graph $H$ as a spanning subgraph, then  $G\in \SL$.
\end{enumerate}
\noindent {\bf Proof}. (a) If $|Y|>|X|$, then since $H$ is a bipartite graph and $d_H(v)\ge 3$ for any $v\in Y$, 
there is at least one vertex (say $x$) in $X$ such that $d_H(x)\ge 4$. Let $N_H(x)=\{y_1 ,y_2, y_3, y_4, \cdots\}$.
Since $H$ is a bipartite graph, $\cup_{i=1}^4N_H(y_i)\subseteq X$. Since $|N_H(y_i)|\ge 3$ ($1\le i\le 4$) and $|X|\le 7$, 
there are at least two vertices (say $y_1$ and $y_2$) in $\{y_1, y_2, y_3, y_4\}$ 
such that $(N_H(y_1)-\{x\})\cap (N_H(y_2)-\{x\})\not =\emptyset$.
Let $x_1$ be a vertex in $(N_H(y_1)-\{x\})\cap (N_H(y_2)-\{x\})$. Then $xy_1x_1y_2x$ is a 4-cycle in $H$ with
 $d_H(x)\ge 4$. Theorem 3.3(a) is proved for this case.

Next, we consider the case $|Y|=|X|$. \\
\indent We may assume $H\not \in \SL$. 
Since $|X|\le 7$, $|V(H)|=|X|+|Y|\le 14$. 

 If  $\delta(H)\le 2$, then similar to the argument  above,    $H$ has a 4-cycle  with the stated properties.  
 We are done if $\delta(H)\le 2$. 
Thus, in the following we assume $\delta(H)\ge 3$.

 If $\kappa'(H)\ge 3$, 
then  by Theorem F(b), either $H\in \SL$, contrary to $H\not \in \SL$, 
or $H\in \{P, P_{14}\}$, contrary to that $H$ is a bipartite graph.
Thus $\kappa'(H)\le 2$.

Let $E_1$ be a minimum edge-cut of $H$ with $|E_1|\le 2$. Let $H_1$ and $H_2$ be the two components of $H-E_1$ 
and $|V(H_1)|\le |V(H_2)|$.
Since $\delta(H)\ge 3$ and  $|V(H)|\le 14$, no matter $|E_1|=1$ or 2, we have 
$\delta(H_1)\ge 2$ with $|D_2(H_1)|\le 2$ and $1<|V(H_1)|\le 7$. 
By Theorem F(a), $H_1$ is not reduced, contrary to that $H$ is reduced. Theorem 3.3(a) is proved. 

\vskip 0.1 in

\noindent (b). If $\delta(H)\le 2$, then similar to the argument above, 
$H$ has a 4-cycle  with a vertex of degree at least 4 in $X$.
We are done for this case.

If $\delta(H)\ge 3$, then let $x_0$ be a vertex in $X$. Let $y_1$, $y_2$ and 
$y_3$ be three distinct vertices in $N(x_0)$. Since $H$ is a connected bipartite graph,
 $\cup_{i=1}^3(N_H(y_i)-\{x_0\})\subseteq X-\{x_0\}$ and so $|\cup_{i=1}^3(N_H(y_i)-\{x_0\})|\le |X|-1=5$.
 Since $d_H(y_i)\ge 3$ $(1\le i\le 3)$, 
$|N_H(y_i)-\{x_0\}|\ge 2$. Thus, $\sum_{i=1}^3|N_H(y_i)-\{x_0\})|\ge 6>5\ge |\cup_{i=1}^3(N_H(y_i)-\{x_0\})|$. 
Hence, there are some $i, j\in \{1,2,3\}$ ($i\not=j$) such that
$(N_H(y_i)-\{x_0\})\cap (N_H(y_j)-\{x_0\})\not =\emptyset$, and so $H$ has a 4-cycle. Theorem 3.3(b) is proved.

\vskip 0.1 in

\noindent (c).  Suppose $G\not \in \SL$. 
Let $n=|V(G)|$. 
If $n\ge 16$, then $\frac{n+5}{3}\ge 7\ge |X|$ and $|Y|\ge 9>|X|$. By Lemma 3.2, $G=H$ and $F(G)=3$.
By Theorem 3.1, $G$ has no 4-cycles. But by (a) above,  $G$ has a 4-cycle, a contradiction. Thus  $G\in \SL$ if $n\ge 16$.

If $n\le 14$, then since $\kappa'(G)\ge 3$ and $G\not \in \SL$, by Theorem F(b), $G\in \{P, P_{14}\}$. However,
$P$ and $P_{14}$ have no  spanning bipartite subgraphs with the stated properties. This is impossible.

If $n=15$, then by Theorem F(c), $G$ has girth at least 5.
Since $|X|\le 7$, $|Y|\ge 8>|X|$. By (a) again, $G$ has a 4-cycle, a contradiction.
Theorem 3.3(c) is proved. $\Box$

\vskip 0.1 in

Using  Theorems 3.1 and 3.3, we prove the following result,  Theorem 3.4,  for the size of 
 maximum matchings in reduced graphs, which is an improvement of a result in \cite{CL1}.\\

Let $q(G)$ denote the number of odd components of $G$.\\
{\bf Theorem H} (Berge \cite{B}, Tutte \cite{Tu}).
 Let $G$ be a graph of order $n$. 
Then $\alpha'(G)=(n-t)/2$, where $t=\max_{S\subset V(G)}\{q(G-S)-|S|\}$. $\Box$

\vskip 0.1 in
\noindent {\bf Theorem 3.4}. Let $G$ be a 3-edge-connected reduced graph of order $n$ and $G\not \in \SL$. 
If $n\le 17$, then $\alpha'(G)\ge (n-1)/2$. \\ 
{\bf Proof}. By Theorem F(b), if $n\le 14$, then $G\in \{P, P_{14}\}$ and so $G$ has a  perfect matching.
We are done for $n\le 14$. Thus, we may assume $n\ge 15$.

Let $t$ be the integer defined in Theorem H. By way of contradiction, suppose  $t\ge 2$. 
Let $S\subset V(G)$ be chosen such that $t=q(G-S)-|S|$. Let $m=q(G-S)$ and let $G_1$, $G_2$, $\cdots$, 
$G_m$ be the odd components of $G-S$.  We may assume that
\begin{eqnarray*}
|V(G_1)|\le |V(G_2)|\le \cdots \le |V(G_m)|.
\end{eqnarray*}

For each odd integer $i$, let ${\cal R}_i$ be the collection of components of $G-S$
 consisting of exactly $i$ vertices, and 
let $r_i=|{\cal R}_i|$.
Let $S_i=\cup_{H\in {\cal R}_i}V(H)$. Then $|S_i|=ir_i$ ($i=1, 3,\cdots$). 
 For each component $H$ of $G-S$, let $\partial (H)$ 
be the set of edges in which every edge incident with at least one vertex in $V(H)$.  
Then 
\begin{eqnarray}
 n&\ge &|S|+\Sigma_{i=1}^m|V(G_i)|=|S|+r_1+3r_3+5r_5+\cdots;\\
 m&=&|S|+t=q(G-S)=r_1+r_3+r_5+\cdots. 
\end{eqnarray}
We have 
\begin{eqnarray}
n&\ge&|S|+(r_1+r_3+r_5+\cdots)+(2r_3+4r_5+\cdots);\nonumber \\
n&\ge &|S|+m+2(r_3+2r_5+\cdots)=2|S|+t+2(r_3+2r_5+\cdots).
\end{eqnarray}

By (7), $t\ge 2$ and $n\le 17$, $2|S|\le 17-t\le 15$ and so $|S|\le 7$. 
Furthermore, if $|S|=7$, then by (7) again,
$2(r_3+2r_5+\cdots)= n-2|S|-t\le 1$ and so $r_i=0$ ($i=3,5, \cdots)$. Thus, 
$V(G)=S\cup S_1$.  Since $n\ge 15$, $|S_1|=r_1=n-|S|\ge 8>|S|$.

Let $H$ be the bipartite graph induced by the edges between $S$ and $S_1$.
 Since $\delta(G)\ge 3$ and each vertex $v$ in $S_1$ is only
adjacent to the vertices in $S$, $d_H(v)\ge 3$ for any $v\in S_1$. Therefore,
 $G$ has  a spanning bipartite subgraph $H$
 with the properties stated in Theorem 3.3.
By Theorem 3.3(c),  $G\in \SL$, a contradiction.

In the following, we assume that $|S|\le 6$.

\vskip 0.1 in

\noindent {\bf Case 1}. $r_1+r_3=0$.\\
\indent Let $i\ge 5$ be the smallest integer such that $r_i\not =0$. Then by (5), $m=|S|+t$ and $t\ge 2$,
\begin{eqnarray*}
n\ge |S|+im\ge |S|+5m=6|S|+5t\ge 6|S|+10.
\end{eqnarray*}
Therefore, since $n\le 17$, $|S|\le \frac{n-10}{6}\le \frac{7}{6}$ and so $|S|=1$ and $i=5$. 

Hence, $|V(G_1)|=5$. Let $H=G[S\cup V(G_1)]$. Since $G$ is reduced, $H$ is reduced. 
Since $|S|=1$ and  $G$ is 3-edge-connected, $H$ is a graph with $|V(H)|=|V(G_1)|+|S|=6$ and $\delta(H)\ge 3$. By Theorem F(a), 
$H$ is not reduced, a contradiction. Case 1 is proved.

\vskip 0.1 in

\noindent {\bf Case 2}. $r_1+r_3\not =0$.\\
\indent  Since $G$ is $K_3$-free and $\delta(G)\ge 3$, 
\begin{eqnarray}
 |\partial(H_0)|\ge 3 \mbox{ for each $H_0\in {\cal R}_1$; and }  |\partial(H_1)|\ge 7 \mbox{  for each $H_1\in {\cal R}_3$.}
\end{eqnarray}

Let $G_0=G[S_0\cup S_1\cup S_3]$ where $S_0$ is the largest subset of $S$ such that $G_0$ is connected.
 Then $|S_0|\le |S|$ and $\displaystyle E(G_0)=\cup_{H\in {\cal R}_1\cup {\cal R}_2}E(G[V(H)\cup S_0])$. 
 By Theorem E(a), $G_0$ is a reduced graph with 
\begin{eqnarray}
 |V(G_0)|=|S_0|+|S_1|+|S_3|=|S_0|+r_1+3r_3. 
\end{eqnarray}
Since for any two $H_1$, $H_2\in {\cal R}_1\cup {\cal R}_3$ with $H_1\not =H_2$, 
$\partial(H_1)\cap \partial (H_2)=\emptyset$, 
 $|E(G_0)|=\Sigma_{H\in {\cal R}_1\cup {\cal R}_2}|\partial(H)|+|E(G[S_0])|$. By (8) 
\begin{eqnarray}
|E(G_0)|\ge 3r_1+7r_3. 
\end{eqnarray}

\noindent {\bf Claim 1}. $G_0\not \in \{K_1, K_2, K_{2,s}\}$ $(s\ge 1)$.\\
\indent Since  each vertex $v\in S_1$ is only adjacent to the vertices in $S$ and
 each vertex $v\in S_3$ is only adjacent to vertices in $S\cup S_3$, and since $\delta(G)\ge 3$, 
$d_H(v)=d(v)\ge 3$ for any $v\in S_1\cup S_3$, and so $|S|\ge 3$. 
Thus $G_0\not \in \{K_1, K_2\}$. Next we will show $G_0\not =K_{2,s}$.

\vskip 0.1 in
 Suppose that $G_0=K_{2,s}$ ($s\ge 1$).
 Then $G_0$ has at most two vertices of degree greater than 2.
Thus $r_3=|S_3|=0$ and $r_1=|S_1|\le 2$. By (5), (6), $t\ge 2$ and $m=|S|+t$,
\begin{eqnarray*}
n\ge |S|+r_1+5(m-r_1)=|S|+5m-4r_1=6|S|+5t-4r_1\ge 6|S|+2.
\end{eqnarray*}
Since $n\le 17$, $6|S|\le n-2\le 15$. Thus, $|S|\le 2$, contrary to $|S|\ge 3$. Claim 1 is proved.

\vskip 0.1 in

Since $G_0\not \in \{K_1, K_2, K_{2,s}\}$ $(s\ge 1)$, by  Theorem E(c), $F(G_0)\ge 3$. By Theorem E(b),
 $|E(G_0)|\le 2|V(G_0)|-5$.  By (9) and (10),
\begin{eqnarray}
3r_1+7r_3&\le &|E(G_0)|\le 2|V(G_0)|-5=2(|S_0|+r_1+3r_3)-5, \nonumber \\
     r_1+r_3&\le& 2|S_0|-5\le 2|S|-5.
\end{eqnarray}
By (5), (6), (11), $n\le 17$ and $t\ge 2$,
\begin{eqnarray*}
  n&\ge &|S|+r_1+3r_3+5(m-r_1-r_3)\ge 6|S|+5t-2(r_1+r_3)-2r_1;\nonumber \\
   2r_1&\ge & 6|S|+5t-2(r_1+r_3)-n\ge 6|S|-2(r_1+r_3)-7.  \nonumber \\
   2r_1&\ge& 6|S|-2(2|S|-5)-7=2|S|+3
\end{eqnarray*}

Therefore, $r_1\ge |S|+2$. By (11) and $|S|\le 6$,
\begin{eqnarray*}
|S|+2+r_3&\le& r_1+r_3\le 2|S|-5; \\
    r_3&\le & |S|-7\le -1,
\end{eqnarray*}
contrary to $r_3\ge 0$.  Theorem  3.4 is proved. $\Box$

\vskip 0.1 in

\noindent {\bf 4. Degree condition of adjacent vertices for supereulerian graphs}
\vskip 0.1 in

With the theorems on the properties of reduced graphs proved in section 3, we are able to prove a new result for 3-edge-connected graph $G$ that 
satisfies $\overline{\sigma}_2(G)\ge \frac{2n}{p}-2$.

\vskip 0.1 in

 Different from the study on Ore-type degree sum conditions of nonadjacent vertices for hamiltonian graphs,  
Brualdi and Shaney \cite{Brualdi} studied  degree-sum  conditions of adjacent vertices 
to obtain a result on Hamiltonian line graphs.

\noindent {\bf Theorem I} (Brualdi \cite{Brualdi}). 
Let $G$ be a graph of order $n\ge 4$. If for any edge $uv\in E(G)$, $\overline{\sigma}_2(G)\ge n$,
then $G$ contains a dominating Eulerian subgraph, hence $L(G)$ is hamiltonian.

\vskip 0.1 in
Since then, many results had been found on the 
degree-sum conditions 
of adjacent vertices for spanning and dominating Eulerian subgraphs of graphs
(see \cite{BCKV86, chen1, CL2, LLZ,V1}).
The following was proved by Veldman  \cite{V1}. 

\noindent
{\bf Theorem J} (Veldman \cite{V1}).
Let $G$ be a 2-edge-connected simple graph of order $n$. 
If for any $uv \in E(G)$, $\overline{\sigma}_2(G)>\frac{2n}{5} - 2,$
then for $n$ sufficiently large, $L(G)$ is Hamiltonian.

\vskip 0.1 in

For 3-edge-connected graphs, the degree-sum condition in Theorem J can be lower.\\
\noindent
{\bf Theorem K} (Chen and Lai \cite{CL2} and Veldman \cite{V1}). 
Let $G$ be a 3-edge-connected simple graph of order $n$. If $n$ is large and 
$\overline{\sigma}_2(G)\ge \frac{n}{5} - 2,$
then  either $G\in \SL$ or $n=10s$ ($s>0$) and $G'=P$ 
with the preimage of each vertex in $P$  is a $K_s$ or $K_s-e$ for some $e\in E(K_s)$.

\vskip 0.1 in

Here we show  how such graphs change
 when $p$ is increased to 15. \\
\noindent {\bf Theorem 4.1}. Let $G$ be a 3-edge-connected simple graph of order $n$. 
If $n$ is sufficiently large and  
\begin{eqnarray}
\overline{\sigma}_2(G)> 2(\frac{n}{15}-1),
\end{eqnarray}
then either $G\in \SL$ or $G' \in\{P, P_{14}\}$. Furthermore, if $\overline{\sigma}_2(G)\ge 2(\frac{n}{14}-1)$ and
$G'=P_{14}$, then $n=14s$ and each vertex in $P_{14}$ is obtained 
by contracting a $K_s$ or $K_s-e$ for some $e\in E(K_s)$.

\vskip 0.1 in
We prove the following lemma first:\\
\noindent {\bf Lemma 4.2}. Let $G$ be a 3-edge-connected graph of order $n$ with $\overline{\sigma}_2(G)\ge \frac{2n}{p}-2$, where 
$p$ is a given positive number. Let $G'$ be the reduction of $G$. Let $v$ be a vertex in $G'$ and $H(v)$ be the preimage of $v$.
Then  when $n$ is large, each of the following holds:\\
(a) If $|V(H(v))|=1$, then for any $x\in N_{G'}(v)$, $|V(H(x))|\ge \overline{\sigma}_2(G)+1-d_{G'}(v)-d_{G'}(x)$.\\
(b) If $|V(H(v))|>1$, then $|V(H(v))|\ge \frac{\overline{\sigma}_2(G)}{2}+1.$\\
\noindent {\bf Proof}. 
For a vertex $y\in V(G)$, let $i(y)$ be the number of edges in $G'$ incident with $y$ in $G$.
 If $y\in V(H(v))$ where $H(v)$ is the preimage of $v\in V(G')$, then 
\begin{eqnarray}
d_G(y)\le i(y)+|V(H(v))|-1.
\end{eqnarray}
By Theorem D, $|V(G')|\le 3p-4$. Then
\begin{eqnarray}
    \Delta(G')\le |V(G')|-1\le 3p-5.
\end{eqnarray}

\vskip 0.1 in
\noindent (a)\  Since $|V(H(v))|=1$, $v$ is a trivial contraction. Then $d_{G'}(v)=d_G(v)$.
For any $x\in N_{G'}(v)$,  there is a vertex $x_0$ in $G$ such that $e=xy=x_0v$.
Then $d_G(x_0)\le d_{G'}(x)+|V(H(x))|-1$. Since $d_G(v)+d_G(x_0)\ge \overline{\sigma}_2(G)$,
\begin{eqnarray*}
 \overline{\sigma}_2(G)\le d_G(v)+d_G(x_0)\le d_{G'}(v)+d_{G'}(x)+|V(H(x))|-1. 
\end{eqnarray*}
Lemma 4.2(a) is proved. 

\vskip 0.1 in

\noindent (b). Since $|V(H(v))|>1$, $E(H(v))\not =\emptyset$. Let $xy$ be an edge in $E(H(v))$. 
There are at most $d_{G'}(v)$ number of edges in $E(G')$ incident with 
$x$ and $y$ and so $i(x)+i(y)\le d_{G'}(v)\le \Delta(G')$.
Since   $d_G(x)+d_G(y)\ge \overline{\sigma}_2(G)\ge \frac{2n}{p}-2$,  by (13) and (14),
\begin{eqnarray}
  \overline{\sigma}_2(G)&\le &d_G(x)+d_G(y) \nonumber \\
            &\le& (i(x)+|V(H(v))|-1)+(i(y)+|V(H(v))|-1);\nonumber \\
 \overline{\sigma}_2(G)&\le& i(x)+i(y)+2|V(H(v))|-2;\\
 \frac{2n}{p}-(3p-5)&\le& \overline{\sigma}_2(G)-(i(x)+i(y))+2 \le |V(H(v))|.\nonumber
\end{eqnarray}
Since $p$ is a fixed, when $n$ is large ( say $n>p(3p-5)$), 
$H(v)$ has an edge $xy$  such that $i(x)=i(y)=0$.
Thus, by (15), $|V(H(v))|\ge \frac{\overline{\sigma}_2(G)}{2}+1$. Lemma 4.2(b) is proved. $\Box$

\vskip 0.1 in

\noindent {\bf Proof of Theorem 4.1}.  Suppose  
that $G\not \in \SL$. Let $G'$ be the reduction of $G$. By Theorem A,  $G'\not \in \SL$.
Since $\kappa'(G)\ge 3$, 
 $\kappa'(G')\ge 3$. By Theorem D with $p=15$,  $\alpha'(G')\le 15/2$ and so $\alpha'(G')\le 7$. 
By Theorem E(d),  $|V(G')|\le 3\alpha'(G')-4=17$. 
Thus, by Theorem 3.4, $\alpha'(G)\ge (|V(G')|-1)/2$ and so $|V(G')|\le 15$. 
If $|V(G')|\le 14$, then by Theorem F(b) and $G'\not \in \SL$, $G'\in \{P, P_{14}\}$. We are done for this case.

Next, we show that $|V(G')|=15$ is impossible.

If $|V(G')|=15$, then by Theorem F(c), $G'$ has girth at least 5 and 
 $V(G')=D_3(G')\cup D_4(G')$ where  $D_4(G')$ is an independent set. 
Hence, for any $xy\in V(G')$, 
\begin{eqnarray}
d_{G'}(x)+d_{G'}(y)\le 7.
\end{eqnarray}
\vskip 0.1 in

Let $Y_0=\{v\in V(G')\ | \  |V(H(v))|=1\}$. Let $X=\cup_{v\in Y_0} N_{G'}(v)$. Let $Z=V(G')-X-Y_0$. 

For each $v\in Z$, $|V(H(v))|>1$. By Lemma 4.2(b) and $\overline{\sigma}_2(G)> 2(\frac{n}{15}-1)$,
\begin{eqnarray}
|V(H(v))|\ge \frac{\overline{\sigma}_2(G)}{2}+1>\frac{n}{15}.
\end{eqnarray}

 For any $x\in X$,  by Lemma 4.2(a), (16) and (12),  $|V(H(x))|\ge \overline{\sigma}_2(G)+1-7>\frac{2n}{15}-8$. 
Since $\cup_{x\in X}V(H(x))\subseteq V(G)$, 
\begin{eqnarray}
n=|V(G)|\ge \Sigma_{x\in X}|V(H(x))|\ge |X|(\frac{2n}{15}-8).
\end{eqnarray}
Thus,  when $n$ is large, $|X|\le 7$.

\vskip 0.1 in

\noindent {\bf Case 1}. $|Z|\le 1$.  Let $Y=Y_0\cup Z$. Note that if $|Z|=1$, then by the definitions of $Z$ 
 the vertex in $Z$ is only adjacent to 
vertices in $X$. Thus, the edges between $X$ and $Y$ forms a spanning bipartite subgraph $H_a$ of $G$
such  that $d_{H_a}(v)=d(v)\ge 3$ for any $v\in Y$.
Since $|X|\le 7$ and $|X|+|Y|=|V(H_a)|=|V(G')|=15$, $|X|<|Y|$.
Thus, $H_a$ is a bipartite graph with the properties stated in Theorem 3.3(a), and so 
$H_a$ has 4-cycle, contrary to that $G'$ has girth at least 5.  
Case 1 is proved. 

\vskip 0.1 in
\noindent {\bf Case 2}. $|Z|\ge 2$. Then $|X|+|Y_0|\le 13$. 
Let $H_b$ be the bipartite subgraph formed by the edges between $X$ and $Y_0$. 
Since $\kappa'(G')\ge 3$, $d_{H_b}(v)=d(v)\ge 3$ for any $v\in Y_0$. 
Since  $V(G)=\cup_{x\in X}V(H(x))\cup Y_0\cup_{v\in Z}V(H(v))$ and  $|Z|=15-(|X|+|Y_0|)$,
by (17) and (18)
\begin{eqnarray}
 n=|V(G)|&\ge& |X|(\frac{2n}{15}-8)+|Y_0|+|Z|\frac{n}{15}\\
         &= & |X|\frac{2n}{15}-8|X|+|Y_0|+(15-|X|-|Y_0|)\frac{n}{15}\nonumber \\
         &\ge & 
          n+\frac{|X|-|Y_0|}{15}-8|X|+|Y_0|.\nonumber
\end{eqnarray}
Therefore, when $n$ is large, $|X|\le |Y_0|$. Since $|X|+|Y_0|\le 13$, $6\ge |X|$.  
$H_b$ is a bipartite graph with the properties stated in Theorem 3.3. 
By Theorem 3.3(a) and (b), $H_b$ has a 4-cycle, a contradiction. 
This shows that $|V(G')|=15$ is impossible.

\vskip 0.1 in

Next, we assume that $\overline{\sigma}_2(G)\ge \frac{2n}{14}-2$ and $G'=P_{14}$.
\vskip 0.1 in

\noindent {\bf Claim 1}. $Y_0=\emptyset$.\\
\indent Suppose  $Y_0\not =\emptyset$. Then $X\not =\emptyset$. 
By Lemma 4.2, for each $v\in Z$, $|V(H(v))|\ge \frac{n}{14}$, and for each $x\in X$, $|V(H(x))|\ge \frac{2n}{14}-7$.
Replacing $\frac{n}{15}$ by $\frac{n}{14}$ and replacing $|X|(\frac{2n}{15}-8)$ by $|X|(\frac{2n}{14}-7)$ and using
$|Z|=14-|X|-|Y_0|$ in (19), we have $|X|\le |Y_0|$ when $n$ is sufficiently large.

 Let $H_b$ be the bipartite subgraph defined in Case 2 above.
Since $d_{H_b}(v)\ge 3$ for any $v\in Y_0$,
$|X|\ge 3$. Since  $|Y_0|\ge |X|$ and $G'$ has no $K_{3,3}$, $|Y_0|\ge 4$. 
By Lemma 4.2(a), $Y_0$ is an independent set.
However, by observation on $P_{14}$, 
$|X|=|\cup_{v\in Y_0} N_{G'}(v)|\ge 7$ for any independent set $Y_0$ with $|Y_0|\ge 4$.   $P_{14}$ has no such bipartite subgraph $H_b$.
 Claim 1 is proved.
\vskip 0.1 in

Therefore, $Z=V(P_{14})$. Then by $|V(H(v))|\ge \frac{n}{14}$ for each $v\in Z$,
\begin{eqnarray}
n=|V(G)|=\Sigma_{v\in Z}|V(H(v))|\ge |Z|\frac{n}{14}=n.
\end{eqnarray}
Thus the equality of (20) holds and so $|V(H(v))|=\frac{n}{14}$ for any $v\in V(P_{14})$. Let $s=|V(H(v))|=\frac{n}{14}$.
Since for any $uv\in E(G)$, $d(u)+d(v)\ge \overline{\sigma}_2(G)\ge \frac{2n}{14}-2$, $H(v)$ is either $K_s$ or $K_s-e$ for some
 $e\in E(K_s)$. (See $G_b$ in Figure 4.1(b) for  such a graph). 
$\Box$

\vskip 0.1 in
\noindent {\bf Remark}: From Theorem 4.1 and Theorem K, we can see that for a 3-edge-connected graph $G$ of order $n$ with 
$\overline{\sigma}_2(G)\ge \frac{2n}{p}-2$, 
the structures of $G$ change when $p$ is increased:\\
\indent (a) if $\overline{\sigma}_2(G)> \frac{2n}{10}-2$ then $G\in \SL$; \\
\indent (b) if $\overline{\sigma}_2(G)\ge \frac{2n}{10}-2$ then $G\in \SL$ or $G=G_a$ as shown in Figure 4.1(a) \\
\indent \ \ \ \ \  where $n=10s$ and each circle is a $K_s$ or a $K_s-e$;  \\
\indent (c) if $\overline{\sigma}_2(G)> \frac{2n}{14}-2$, then $G\in \SL$ or $G'=P$;  \\ 
\indent (d) if $\overline{\sigma}_2(G)\ge \frac{2n}{14}-2$, then $G\in \SL$ or $G'=P$ or $G=G_b$ as shown in\\
\indent \ \ \ \ \  Figure 4.1(b) where $n=14s$ and each circle is a $K_s$ or $K_s-e$;\\
\indent (e) if  $\overline{\sigma}_2(G)> \frac{2n}{15}-2$, then $G\in \SL$ or $G'\in \{P, P_{14}\}$.

\vskip 0.1 in
Graphs $G_a$ and $G_b$ in Figure 4.1 are the extremal graphs with the boundary value on $p=10$ and 14 for  
$\overline{\sigma}_2(G)\ge \frac{2n}{p}-2$, while $G_c$ is the next possible extremal graph for $p=16$.

$$
\setlength{\unitlength}{.80pt}
\begin{picture}(400,80)(0,-5)
\put(-5,70){$G_a$}     
\multiput(20,50)(14,0){2}{\circle*{4}}
\multiput(20,50)(14,0){2}{\circle{9}}

\put(86,50){\circle*{4}}
\put(86,50){\circle{9}}
\put(101,50){\circle*{4}}
\put(101,50){\circle{9}}

\multiput(20,50)(14,0){2}{\line(1,0){20}}
\put(60,80){\circle*{4}}
\put(60,80){\circle{9}}

\put(34,15){\circle*{4}}
\put(34,15){\circle{9}}
\multiput(46,30)(28,0){2}{\circle*{4}}
\multiput(46,30)(28,0){2}{\circle{9}}

\put(86,15){\line(2,5){14}}
\put(86,50){\line(1,0){15}}
\put(20,50){\line(2,-5){14}}
\put(20,50){\line(4,3){40}}
\put(60,65){\circle*{4}}
\put(60,65){\circle{9}}

\put(60,65){\line(0,1){15}}
\put(46,30){\line(2,5){14}}
\put(60,65){\line(2,-5){14}}
\put(46,30){\line(2,1){40}}
\put(34,50){\line(2,-1){40}}
\put(34,50){\line(1,0){52}}

\put(86,15){\circle*{4}}
\put(86,15){\circle{9}}

\put(74,30){\line(4,-5){12}}
\put(34,15){\line(3,4){11}}
\put(34,15){\line(1,0){52}}
 
\put(60,80){\line(4,-3){40}}   


\put(130,70){$G_b$}     
\multiput(150,50)(14,0){2}{\circle*{4}}  
\multiput(150,50)(14,0){2}{\circle{9}}

\put(216,50){\circle*{4}}
\put(216,50){\circle{9}}
\put(231,50){\circle*{4}}
\put(231,50){\circle{9}}

\multiput(150,50)(14,0){2}{\line(1,0){20}} 
\put(190,80){\circle*{4}}
\put(190,80){\circle{9}}
\put(164,15){\circle*{4}}
\put(164,15){\circle{9}}

\multiput(176,30)(28,0){2}{\circle*{4}}
\multiput(176,30)(28,0){2}{\circle{9}}
\put(231,30){\circle*{4}}
\put(231,30){\circle{9}}
\put(231,70){\circle*4}
\put(231,70){\circle{9}}

\put(251,40){\circle*4}
\put(251,40){\circle{9}}

\put(251,60){\circle*4}
\put(251,60){\circle9}

\put(251,40){\line(-2,-1){20}}
\put(251,40){\line(-2,1){20}}
\put(251,40){\line(-2,3){20}}

\put(251,60){\line(-2,-1){20}}
\put(251,60){\line(-2,1){20}}
\put(251,60){\line(-2,-3){20}}

\put(216,15){\line(1,1){15}}
\put(216,50){\line(1,0){15}}
\put(150,50){\line(2,-5){14}}
\put(150,50){\line(4,3){40}}
\put(190,65){\circle*{4}}
\put(190,65){\circle9}

\put(190,65){\line(0,1){15}}
\put(176,30){\line(2,5){14}}
\put(190,65){\line(2,-5){14}}
\put(176,30){\line(2,1){40}}
\put(164,50){\line(2,-1){40}}
\put(164,50){\line(1,0){52}}

\put(216,15){\circle*{4}}
\put(216,15){\circle9}

\put(204,30){\line(4,-5){12}}
\put(164,15){\line(3,4){11}}
\put(164,15){\line(1,0){52}}
 
\put(190,80){\line(4,-1){40}}


\setlength{\unitlength}{1pt}

\put(230,40){\circle*3}
\put(230,40){\circle7}

\put(230,40){\line(1,-3){10}}
\put(230,40){\line(1,0){10}}
\put(230,40){\line(3,2){30}}    

\put(240,40){\circle*3}
\put(240,40){\circle7}

\put(240,40){\line(1,0){40}}
\put(240,40){\line(3,-2){30}}

\put(280,40){\circle*3}
\put(280,40){\circle7}

\put(240,10){\circle*3}
\put(240,10){\circle7}

\put(240,10){\line(1,1){10}}
\put(240,10){\line(1,0){40}}

\put(280,10){\circle*3}
\put(280,10){\circle7}

\put(280,10){\line(-1,1){10}}
\put(280,10){\line(1,0){20}}

\put(250,20){\circle*3}
\put(250,20){\circle7}

\put(250,20){\line(1,3){10}}
\put(250,20){\line(3,2){30}}

\put(270,20){\circle*3}
\put(270,20){\circle7}

\put(260,50){\circle*3}
\put(260,50){\circle7}

\put(260,50){\line(0,1){10}}
\put(260,50){\line(1,-3){10}}

\put(260,60){\circle*3}
\put(260,60){\circle7}

\put(300,10){\circle*3}   
\put(300,10){\circle7} 

\put(300,10){\line(0,1){20}}
\put(300,10){\line(2,3){20}}

\put(300,30){\circle*3}
\put(300,30){\circle7}

\put(290,40){\circle*3}
\put(290,40){\circle7}

\put(290,40){\line(-1,0){10}}
\put(290,40){\line(1,1){10}}
\put(290,40){\line(1,-1){10}}

\put(300,50){\circle*3}
\put(300,50){\circle7}

\put(300,60){\circle*3}
\put(300,60){\circle7}

\put(300,60){\line(-1,0){40}}
\put(300,60){\line(0,-1){10}}
\put(300,60){\line(1,-1){20}}

\put(310,40){\circle*3}
\put(310,40){\circle7}

\put(310,40){\line(-1,1){10}}
\put(310,40){\line(-1,-1){10}}
\put(310,40){\line(1,0){10}}

\put(320,40){\circle*3}
\put(320,40){\circle7}

\put(220,57){ $G_c$}

\put(5, -5){(a) \footnotesize $\overline\sigma(G_a)=\frac{2n}{10}-2$ }

\put(120,-5){(b) \footnotesize $\overline\sigma(G_b)=\frac{2n}{14}-2$}

\put(240, -5){(c) \footnotesize $\overline\sigma(G)=\frac{2n}{16}-2$}
\put(140, -22) {\small Figure 4.1}
\end{picture}
$$

\vskip 0.1 in

\vskip 0.1 in
Let $G$ be the graph of order $n$ defined in Figure 4.1(c) in which  each circle is a $K_{n/16}$. 
Then $\overline\sigma_2(G)\ge 2(\frac{n}{16}-1)$ and $G'=P_{16}$. 
Thus, (12) in Theorem 4.1 cannot be replaced by  $\overline\sigma_2(G)\ge\frac{2n}{16}-2$. But
if  Conjecture B is true, then we can have $\overline\sigma_2(G)\ge \frac{2n}{16}-2$ for (12) with the conclusion that
 either $G\in \SL$ or $G'\in \{P, P_{14}, P_{16}\}$ and when $G'=P_{16}$, $G=G_c$.

\vskip 0.1 in

As we can see that $P_{14}$ and $P_{16}$ can be contracted to $P$ by contracting a subgraph to a vertex in $P$.
If we relax the conclusion of Theorem 4.1 from ``the reduction of $G$ is in $\{P, P_{14}\}$'' to 
``$G$ can be contracted to $P$", 
the degree condition (12) may be lower.
It was conjectured in \cite{CL1} that for any 3-edge-connected graph $G$ 
of order $n$ if $\overline\sigma_2(G)> n/9-2$, 
then when $n$ is large either $G\in \SL$ or $G$ can be contracted to $P$.

\end{document}